\documentclass[conference]{IEEEtran}
\IEEEoverridecommandlockouts
\usepackage{cite}
\usepackage{amsmath,amssymb,amsfonts}
\usepackage{algorithmic}
\usepackage{graphicx}
\usepackage{textcomp}
\usepackage{xcolor}
\usepackage{hyperref}

\usepackage{subfig}

\def\BibTeX{{\rm B\kern-.05em{\sc i\kern-.025em b}\kern-.08em
    T\kern-.1667em\lower.7ex\hbox{E}\kern-.125emX}}
    
\definecolor{babyblue}{rgb}{0.54, 0.81, 0.94}
    
\begin{document}

\title{Solving Sensor Placement Problems In Real Water Distribution Networks Using Adiabatic Quantum Computation
}

\author{\IEEEauthorblockN{1\textsuperscript{st} Stefano Speziali}
\IEEEauthorblockA{\textit{Idea-re S.r.l. } \\
Perugia, Italy \\
sspeziali@idea-re.eu}
\and
\IEEEauthorblockN{2\textsuperscript{nd} Federico Bianchi}
\IEEEauthorblockA{\textit{Idea-re S.r.l. } \\
Perugia, Italy \\
fbianchi@idea-re.eu}
\and
\IEEEauthorblockN{3\textsuperscript{rd} Andrea Marini}
\IEEEauthorblockA{\textit{Idea-re S.r.l. } \\
Perugia, Italy \\
amarini@idea-re.eu}
\and
\IEEEauthorblockN{4\textsuperscript{th} Lorenzo Menculini}
\IEEEauthorblockA{\textit{Idea-re S.r.l. } \\
Perugia, Italy \\
lmenculini@idea-re.eu}
\and
\IEEEauthorblockN{5\textsuperscript{th} Massimiliano Proietti}
\IEEEauthorblockA{\textit{Idea-re S.r.l. } \\
Perugia, Italy  \\
mproietti@idea-re.eu}
\and
\IEEEauthorblockN{6\textsuperscript{th} Loris F.  Termite}
\IEEEauthorblockA{\textit{Idea-re S.r.l.} \\
Perugia, Italy \\
ltermite@idea-re.eu}
\and
\IEEEauthorblockN{7\textsuperscript{th} Alberto Garinei}
\IEEEauthorblockA{\textit{Department of Engineering Sciences,} \\
\textit{ Guglielmo Marconi University}\\
\textit{\& Idea-re S.r.l. }\\
Perugia, Italy \\
a.garinei@unimarconi.it}
\and
\IEEEauthorblockN{8\textsuperscript{th} Marcello Marconi}
\IEEEauthorblockA{\textit{Department of Engineering Sciences,} \\
\textit{ Guglielmo Marconi University}\\
\textit{\& Idea-re S.r.l. }\\
Perugia, Italy \\
m.marconi@unimarconi.it}
\and
\IEEEauthorblockN{9\textsuperscript{th} Andrea Delogu}
\IEEEauthorblockA{\textit{BlueGold S.r.l.} \\
Milan, Italy \\
andrea@blue-gold.it}
}

\maketitle

\thispagestyle{plain}
\pagestyle{plain}

\begin{abstract}
Quantum annealing has emerged in the last few years as a promising quantum computing approach to solving large-scale combinatorial optimization problems. In this paper, we formulate the problem of correctly placing pressure sensors on a Water Distribution Network (WDN) as a combinatorial optimization problem in the form of a Quadratic Unconstrained Binary Optimization (QUBO) or Ising model. Optimal sensor placement is indeed key to detect and isolate fault events. We outline the QUBO and Ising formulations for the sensor placement problem starting from the network topology and few other features. We present a detailed procedure to solve the problem by minimizing its Hamiltonian using PyQUBO, an open-source Python Library. We then apply our methods to the case of a real Water Distribution Network. Both simulated annealing and a hybrid quantum-classical approach on a D-Wave machine are employed. 
\end{abstract}


\section{Introduction}

Water leaks in Water Distribution Networks (WDNs) can be cause of significant economic loss,  besides being waste of important resources \cite{farley2003losses, adedeji2017towards}.  Asset management of WDNs is indeed a relevant issue for the scientific community and novel and more efficient solutions to detect and isolate leaks are always needed.  In particular,  the optimal placement of sensors is crucial if we want to monitor the behavior of a WDN and prevent fault events. 

It is often the case that the number of nodes in a WDN is much larger than the number of available sensors. While the former can be in the order of thousands for a realistic WDN,  there is usually just a few dozens of the latter.  Hence, sensors must be placed such that network-wide global relevant information can be provided.

Several methods to deal with leak detection and isolation have been developed in recent years. See \cite{puust2010review} and references therein for a comprehensive overview.  Standard methods to detect and isolate water leaks employ the use of genetic algorithms \cite{perez2014leak, casillas2013optimal}, while more modern techniques contemplate the use of Artificial Intelligence-inspired methods.  For example, in \cite{irofti2017dictionary, irofti2020fault} Dictionary Learning strategies are used to handle pipe leakage in WDNs through fault detection and isolation mechanisms, while in \cite{javadiha2019leak} Convolutional Neural Networks are used to learn different pressure maps characterizing leaks localization.

In this paper\footnote{The study presented here is part of the project ``WATER A.I. - empower the efficiency of WATER 	networks through Artificial Intelligence for IoT'' financed to BlueGold S.r.l by BANDO INNODRIVER S3 Azione I.1.b.1.1– sostegno all’acquisto di servizi per l’innovazione tecnologica, strategica, organizzativa e commerciale delle imprese, ID Progetto 1734422 | CUP E47B20000590007.},  we take a different route and formulate the problem of sensor placement as a combinatorial optimization problem \cite{CombinatorialOptimization}.  Roughly speaking, combinatorial optimization refers to the computation of maxima or minima of a function over a discrete domain. 

Some combinatorial problems can easily become intractable when the number of variables is large enough.  However, it turns out that many of these problems can be addressed by means of a new computational technique, known as Ising machine.  As a matter of fact, in 2011 D-Wave Systems announced D-Wave One,  the first commercially available quantum annealer, able to solve combinatorial optimization problems \cite{johnson2011quantum}. Since then, large-scale combinatorial problems in daily life have been addressed by means of Quantum Annealing (QA) on D-Wave quantum or hybrid quantum-classical solvers with applications ranging from delivery \cite{borowski2020new}, to traffic flow \cite{neukart2017traffic} or job scheduling \cite{venturelli2015quantum}.

The reasons for thinking of the sensor placement as a combinatorial problem are diverse.  First of all, we do not need to simulate the hydraulic behavior of the water network, the mathematical formulation of which can be sometimes very challenging and poorly accurate.  Second,  as stated above, combinatorial problems can be run on Ising machines (like the D-Wave quantum annealer) and often good solutions can be found quickly. This can be useful when we want to monitor large WDN's and optimal sensor configurations are needed in a relatively short time. 

In order to employ Ising machines to solve an optimization problem, one should define the energy function (Hamiltonian) of the Ising model or QUBO (Quadratic Unconstrained Binary Optimization) problem which corresponds to the function we want to minimize (or maximize) \cite{lucas2014ising}. 

We propose in this paper an instance of Hamiltonian whose ground state encodes the optimal sensor placement for a generic WDN.  The lowest energy state is then found by means of Simulated Annealing (SA) or Adiabatic Quantum Computation (AQC). In order to program the optimization, we used PyQUBO, an open source Python library \cite{zaman2021pyqubo} useful to construct QUBOs from the objective functions and constraints of optimization problems.  As we shall see, the algorithm we propose is very flexible and suitable for refinements or generalizations. 

The paper is organized as follows.  In section \ref{Combinatorial Optimization} we review some basic facts about combinatorial problems and discuss how they can be formulated as optimization problems (QUBO).  In section \ref{From QUBO to Ising models and Quantum Adiabatic Optimization} we review how QUBO problems can be mapped to Ising models and briefly discuss quantum adiabatic optimization.  Sections \ref{QUBO formulation of the sensor placement problem} and \ref{Case study: simulations and results on a real WDN} are the core of the paper.  We formulate the problem of how to place sensors on a generic WDN as a QUBO problem and we perform a simulation on a real WDN.  We give our conclusions in section \ref{Conclusion}.  In Appendix \ref{Annealing Methods and Adiabatic Evolution} we give an overview on the basics of annealing methods in Quantum Mechanics and quantum adiabatic evolution, while in Appendix \ref{A technical problem and its resolution} we discuss a possible technical problem that might arise when physically installing sensors and its resolution.

\section{Combinatorial Optimization}\label{Combinatorial Optimization}

As stated in the introduction,  many combinatorial optimization problems can be solved using quantum algorithms.  Before discussing how this can be done, let us review a few basic facts about combinatorial optimization.

A combinatorial optimization problem can be formulated as follows.  Let
\begin{equation}
\mathcal{C}: \mathcal{S} \rightarrow \mathbb{R}
\end{equation}
be a function over a set $\mathcal{S}$ of decision variables to the field of real numbers. We seek to find a (possibly unique) element of $\mathcal{S}$,  $x^* \in \mathcal{S}$,  which minimizes $\mathcal{C}$:
\begin{equation}\label{combinatorial optimization}
x^* = \text{arg min}_{x} \mathcal{C}(x) \, ,
\end{equation}
such that a set of constraints of the form
\begin{equation}\label{generic constraints}
g_{l}(x) = 0 \, , \quad h_{m}(x) \leq 0
\end{equation}
is satisfied.  Here $l = 1, \dots, L$ and $m = 1, \dots, M$ are indices counting the number of equality and inequality constraints, respectively.

The function $\mathcal{C}$ is usually referred to as the ``cost function'' of the problem.  Remarkably, equation \eqref{combinatorial optimization} can be re-written as an unconstrained optimization problem where constraints, as in \eqref{generic constraints},  are given in the form of penalties.  More in detail, we want to find $x^*$ such that
\begin{equation}
x^* = \text{arg min}_{x} H_P(x) \, ,
\end{equation}
with
\begin{equation}
H_P = \mathcal{C}(x) + \sum_l a_l g_l(x)^2 + \sum_m b_m \text{max} [h_m(x),  0] \, .
\end{equation}
Here $a_l$ and $b_m$ are positive real numbers enforcing the constraints $g_l(x) = 0$ and $h_{m}(x) \leq 0$.  The coefficients $a_l$ and $b_m$ must be tuned so that it is disadvantageous to violate any of the constraints. 

When the set of decision variables, $\mathcal{S}$,  is given in terms of binary variables, $\mathcal{S} = \mathbb{B}^n$ with $\mathbb{B} = \{0, 1\}$, and the ``energy function'' $H_P$ is quadratic\footnote{If the problem of interest has higher order interactions, it is still possible to reduce it to a QUBO by introducing ancillary variables. For example, a local expression of the form $x_1 x_2 x_3$ is reduced to $x_1 x_4$ if we define $x_4 = x_2 x_3$.  The latter condition can be imposed by adding to the cost function the term $$3 x_4 + x_2 x_3 - 2 x_2 x_4 -2x_3 x_4.$$ This expression gives a vanishing contribution only when $x_4 = x_2 x_3$ and is positive (penalty) otherwise.},  the problem is usually referred to as QUBO.

QUBO problems are closely related and equivalent to Ising models, where the decision variables are spin, taking values in $\{-1, 1 \}$ instead of $\mathbb{B}$. To go from QUBO to Ising we just need to perform a straightforward change of a variables.  More details are given in the following section.

Thus,  solving a combinatorial optimization problem can be rendered, in many cases, equivalent to finding the ground state of a quantum system, where physics methods can be readily applied.  One of the most promising methods to find the ground state of a given Hamiltonian is provided by AQC, a framework of quantum computation that relies on the adiabatic condition (slow change of parameters), and such that the quantum system is supposed to follow in the time evolution the instantaneous ground state. 

AQC is often used in the literature as synonym of Quantum Annealing (QA).  In fact, the latter defines a broader concept, where also non-adiabatic changes of the underlying Hamiltonian are possible. However,  in this paper we will restrict ourselves to the case where the two concepts can be used interchangeably.

Let us now see how QUBO problems can be mapped to Ising models and discuss briefly how AQC is useful to our purposes. 

\section{From QUBO to Ising models and Quantum Adiabatic Optimization}\label{From QUBO to Ising models and Quantum Adiabatic Optimization}

A classical QUBO problem is usually written as a quadratic function of a set of $n$ binary variables $x_i$ which take value 0 or 1:
\begin{equation}\label{classical QUBO}
H_P = \sum_{i, j} Q_{ij} x_i x_j + \sum_{i} c_i x_i \, .
\end{equation}
Here $Q_{ij}$ and $c_i$ are real numbers. Of course, $Q_{ij} = Q_{ji}$, as the term $x_i x_j$ is symmetric with respect to the exchange of $i$ and $j$. 

The model \eqref{classical QUBO} can be re-written in the form of a Ising model if we perform the straightforward change of variables $x_i = (s_i +1)/2$.  Here, $s_i$ are spin variables taking values $\pm 1$. We get for $H_P$:
\begin{equation}\label{classical ising}
H_P = - \sum_{i, j} J_{ij} s_i s_j - \sum_{i} h_i s_i \, ,
\end{equation}
where $h_i$ and $J_{ij}$ are real numbers related to $c_i$ and $Q_{ij}$ by the following expressions
\begin{equation}
h_i = - \frac{1}{2} \left( c_i + \sum_j Q_{i j} \right) \, , \qquad J_{ij} = - \frac{1}{4}Q_{ij} \, .
\end{equation}
A constant energy-shift in \eqref{classical ising} has been neglected: In general, constant shifts of the energy function can be ignored as they play no role in the minimization (or maximization) procedure.

The quantum version of the same Hamiltonian $H_P$ is simply obtained by replacing $s_i \rightarrow \sigma_{i}^z$, 
\begin{equation}\label{quantum ising}
H_P = - \sum_{i < j} J_{ij} \sigma_i^z \sigma_j^z - \sum_{i} h_i \sigma_i^z \, ,
\end{equation}
where $\sigma_{i}^z = \mathbb{I}^{\otimes(i-1)} \otimes \begin{pmatrix}
1 & 0 \\
0 & -1
\end{pmatrix} \otimes \mathbb{I}^{\otimes(n-i)}$ acts on the $i$-th spin, and with $\mathbb{I}$ being the $2 \times 2$ identity matrix. The Hamiltonian in eqn. \eqref{quantum ising} is easily diagonalized in the basis $\{ |\mathbf{s} \rangle \in \mathbb{C}^{2^n} | \mathbf{s} \in \{ -1,  1\}^{n}  \}$, with $\sigma^z | \pm 1\rangle = \pm | \pm 1 \rangle$. 
When $J_{ij}$ are chosen from a random distribution, the model \eqref{quantum ising} is also referred to as a spin glass,  a well-known NP-hard problem for classical computers.

Finding the ground state of \eqref{classical QUBO} or  \eqref{classical ising} can become a challenging task when the number of variables is large, and here is where annealing methods come to play a role.  We review some basics of annealing methods and the quantum adiabatic theorem in Appendix \ref{Annealing Methods and Adiabatic Evolution}. For more in-depth reviews we refer the reader to e.g. \cite{nakahara2013lectures, tanaka2014physics}.

For the time being, we just point out that an important ingredient for QA is adiabaticity.  According to the quantum adiabatic theorem \cite{messiah1999quantum, amin2009consistency} (see also Appendix \ref{Annealing Methods and Adiabatic Evolution}), the ground state of the model of interest,  defined in this case by some Ising model $H_P$,  can be found by initializing the system in the ground state of some other Hamiltonian $H_B$,  which should be easy to prepare both theoretically and experimentally.  $H_B$ must be chosen such that it does not commute\footnote{This is due to the fact that we want to avoid level crossing: one has to tune (in general) three parameters of a $2 \times 2$ Hamiltonian to get two states with the same energy. We have only one ($t$) and thus we do not expect any degeneracy over the time evolution.} with $H_P$. Then, we perform an adiabatic evolution from $H_B$ to $H_P$,  $H(u) = (1 - u) H_B + u H_P$, where $u = t/T$ is a ``time parameter'' for total evolution time $T$.  In this way we are assured that the solution to the optimization problem encoded in the Hamiltonian $H_P$ is reached after a time $T$, which should be long enough \cite{farhi2000quantum}.

It is customary \cite{farhi2000quantum} to choose $H_B$ to be
\begin{equation}
H_B = - h_1 \sum_i \sigma_i^x \, ,
\end{equation}
where $\sigma_i^x = \mathbb{I}^{\otimes(i-1)} \otimes \begin{pmatrix}
0 & 1 \\
1 & 0
\end{pmatrix} \otimes \mathbb{I}^{\otimes(n-i)}$ again acts on the $i$-th spin, while $h_1$ is a positive real constant.  In this way, the ground state of $H_B$ is an equal superposition of all possible states of $H_P$.  This, in particular, means that no level-crossing, and thus no degeneracy, is expected in the time evolution. 

Let us now move on to formulating the placement of sensors as a QUBO (or Ising) problem.

\section{QUBO formulation of the sensor placement problem}\label{QUBO formulation of the sensor placement problem}

A WDN can be formally interpreted as a graph $G = (\mathcal{V}, E)$, where the nodes (or vertices), denoted collectively as $\mathcal{V}$,  are either tanks or junctions (the former are source of water while the latter distribute the existing water flow to users through the pipes),  whereas the edges, denoted collectively $E$, are the pipes connecting nodes.  Even though we will not need much from graph theory, we refer the reader to \cite{clark1991first} for an overview on graphs and to \cite{vazirani2013approximation} for a nice review on approximation algorithms defined on graph structures.

\subsection{Formulation of the Problem}

Let $x_i$ be a binary variable associated with the $i$-th vertex in $\mathcal{V}$.  We choose $x_i$ to be 1 if the node of our WDN associated with the $i$-th vertex hosts a sensor,  0 otherwise.

The problem of covering our WDN with pressure sensors can be formulated as the minimum vertex cover problem \cite{lucas2014ising}, very well-known in combinatorial optimization. The constraint that every edge\footnote{In the following, an edge will be identified either by a Latin letter,  say $e \in E$, or by its endpoints (say $i$, $j \in \mathcal{V}$) $(ij) \in E$.} $(ij) \in E$ of the graph has at least one vertex associated with a sensor can be encoded in the following Hamiltonian:
\begin{equation}\label{provisional hamiltonian WDN topology}
H_{(0)}' = A \sum_{(ij) \in E} (1 - x_i)(1 - x_j) \,  .
\end{equation}
Then,  we want to minimize the number of vertices with an associated sensor.  In the most general setup,  each node $i \in \mathcal{V}$ comes with an associated cost,  $c_i$ say,  which we always assume to be non-negative \cite{vazirani2013approximation} ($c_i \ge 0$).  Thus,  we want to minimize
\begin{equation}
H_{(1)} = \sum_{i \in \mathcal{V}} c_i x_i \, .
\end{equation}

Given these considerations,  it seems that covering a WDN with pressure sensors amounts to finding the ground state of the hamiltonian $H = H_{(0)}'  + H_{(1)} $.  Of course,  in order to have a well defined minimization problem,  we should specify what the costs $c_i$ correspond to.  We will propose a realistic cost $c_i$ for each node $i \in \mathcal{V}$ below. 

Before defining what $c_i$ will be for us,  it should be kept in mind that in realistic situations the number of sensors is fixed and much smaller than the number of nodes. The constraint that forces the total number of sensors to be equal to a predetermined number $s \in \mathbb{N}$,  i.e.  $\sum_{i \in \mathcal{V}} x_i = s$,  can be included in the minimization problem by considering the following hamiltonian
\begin{equation}\label{constrain fixed number of sensors}
H_{(2)} = B \left( \sum_{i \in \mathcal{V}} x_i - s \right)^2 \, ,
\end{equation}
where $B$ is a positive number, chosen so that the constraint on the fixed number of sensors is always satisfied.

We might also wonder what changes if the number of sensors was not to be strictly equal to a predetermined number, $s$, but rather not greater than that number.  In the latter case the constraint would be instead: $\sum_{i \in \mathcal{V}} x_i \leq s$. 

As we will show in a moment, such a constraint can still be re-written as an equality constraint. However,  an expansion in the number of spins is necessary \cite{lucas2014ising}. 

Let $y_{\alpha}$ for $0 \leq \alpha \leq s$ denote a binary variable, which is 1 if the final number of sensors is $\alpha$ and 0 otherwise.  Then the inequality constraint can be implemented by the following Hamiltonian
\begin{equation}\label{constraint number of sensors inequality}
\tilde{H}_{(2)} = B \left( 1 - \sum_{\alpha=1}^{s} y_{\alpha} \right)^2 + B \left( \sum_{\alpha=1}^s \alpha y_{\alpha} - \sum_{i \in \mathcal{V}} x_i  \right)^2 \, ,
\end{equation}
where the first term is there to ensure that only one among the $y_{\alpha}$'s is equal to 1 (one-hot),  say $y_{\widehat{\alpha}}$, while the second term forces the number of sensors to be equal to $\widehat{\alpha} (< s)$.

One can employ either $H_2$ or $\tilde{H}_2$ for the optimization,  according to the type of constraint considered.

Now,  the fact that the number of sensors is usually much smaller than the number of nodes induces a further complication.  More in detail, not all edges of our network can have at least one sensor associated with one of their endpoint.  Therefore,  to render the problem more consistent,  we should weigh different edges differently, according to their intrinsic properties in the WDN.  We propose to modify the hamiltonian \eqref{provisional hamiltonian WDN topology} in the following fashion:
\begin{equation}
H_{(0)} = A \sum_{(ij) \in E} w_{ij} (1 - x_i)(1 - x_j) \, , 
\end{equation}
where $w_{ij}$ is a (positive) weight for the edge $(ij) \in E$,  which we define later below.

To sum up,  the optimal sensor placement can be formulated as a QUBO problem where the defining hamiltonian is given by
\begin{equation}\label{final hamiltonian}
H_P = H_{(0)} + H_{(1)} + H_{(2)} \, (\text{or } \tilde{H}_2) \, ,
\end{equation}
with $H_{(0)}$,  $H_{(1)}$,  $H_{(2)}$ (or $\tilde{H}_2$) defined as above. 

As a final comment, we point out that one of the four parameters $A$, $B$, $C$ and $D$ could be scaled away: minimizing $H_P$ or $a \times H_P$, with $a$ an arbitrary constant, makes no difference and this gives us the freedom to get rid of one of the four parameters.  However, it is useful to keep all of them so that numerical simulations can be run by choosing the four hyper-parameters independently.

\subsection{Vertex cost function}

We shall now describe what the costs $c_i$,  introduced above,  will be for us.  As remarked before,  the aim of this paper is to derive a simple and realistic model for WDN sensor placement regardless of the hydraulic behavior of the WDN itself.  Two simple and important parameters one may want to take into account are given by the degree of accessibility of the network (some nodes may be physically more accessible than other as they are located in correspondence of control or check valves, fire hydrants and so on) and water consumption in correspondence of the given nodes. 

The former is a useful parameter to consider when physically installing sensors, while the latter guarantees that we pick in the optimization nodes with higher demands of water.

Water consumption is usually given by WDN owners in terms of the amount of water (in a given amount of time) $v_i$ necessary to cover specific needs for users in correspondence of the $i$-th node.  For the following analysis, it is useful to define the dimensionless quantities $\widehat{v}_i = v_i/v_{max}$. Note that $\widehat{v}_i$ lies in the range 0 to 1.

Thus,  we propose to define $c_i$ as
\begin{equation}
c_i = C f_i + D g_i \, ,
\end{equation}
where $f_i \equiv f(\widehat{v}_i)$ is a function of the water need at each node $i$, while $g_i$ weights different nodes differently, according to their degree of accessibility.  $C$ and $D$ are positive weights that can be tuned arbitrarily, but that should not exceed $B$, in order not to violate the constraint \eqref{constrain fixed number of sensors} or \eqref{constraint number of sensors inequality}.  A simple example of accessibility weighting is given by considering $g_i$ to be equal to 1 for nodes that are considered less accessible, and 0 otherwise.  Given also that one may want to monitor nodes with higher rates of consumption,  $f(\widehat{v}_i)$ can be chosen to be a monotonically decreasing function of its argument (normalised between 0 and 1).  Simple models for $f$ might be $f(x) = 1 - x$ or $f(x) = e^{-x}$. 

In section \ref{Case study: simulations and results on a real WDN} we will show the results for sensor placement optimization in the case of simple choices of $f$ and $g$.  Let us now discuss the edge weights $w_{ij}$.

\begin{figure*}[!t]
\centering
\subfloat[Water Distribution Network used for the optimization. It comprises 1368 nodes (blue points) and 1391 edges (solid lines connecting nodes). The solid red point corresponds to the actual tank of the network, while the green point corresponds to a fictitious outlying tank added so to get a more interesting edge betweennes.]{\includegraphics[width=8.4cm]{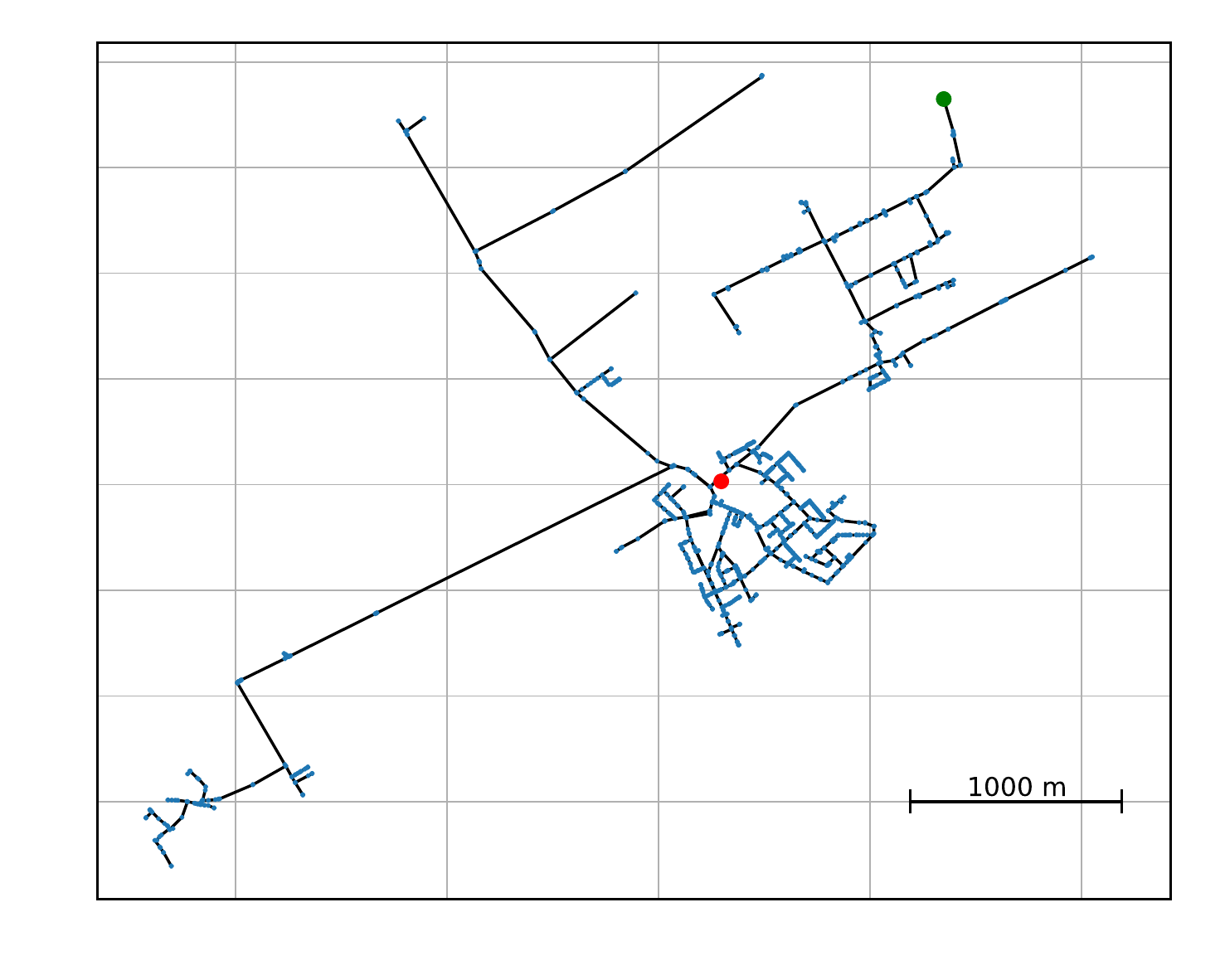}
\label{wdn}}
\hfil
\subfloat[Water Distribution Network tailored centrality. The edge betweennes centrality has been weighted using pipe lengths.  We have normalized it to lie in the range 0 to 1. The darkest edges are those in proximity of the two tanks and where water is more likely to pass through to reach demand nodes of the network.]{\includegraphics[width=8.9cm]{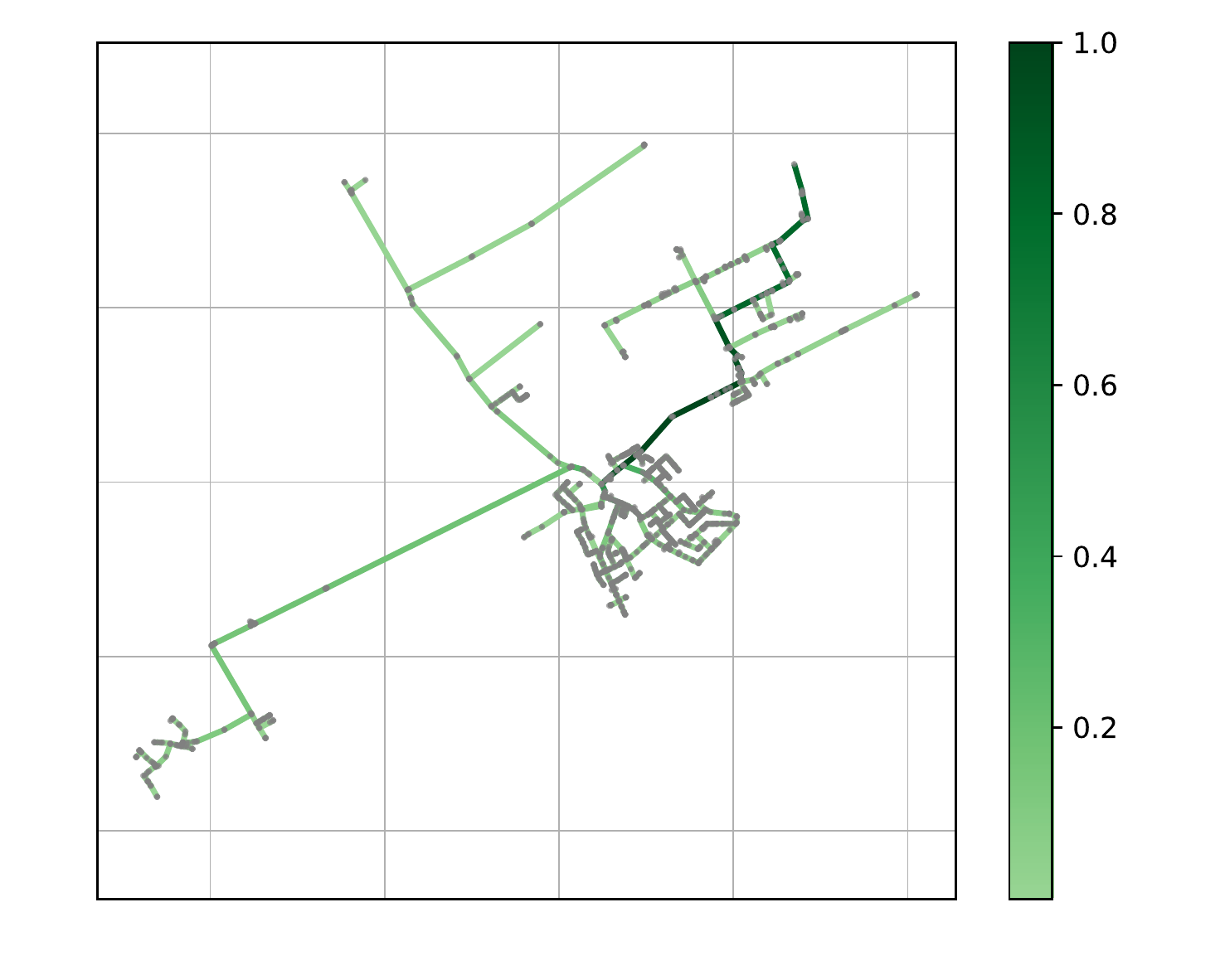}
\label{wdn centrality}}
\footnotesize{\caption{}}
\label{wdn and centrality}
\end{figure*}

\subsection{WDN centrality metrics}

We shall now discuss what $w_{ij}$,  weight for a given edge $(ij) \in E$,  will be for us.  A simple argument suggests that the most ``central'' nodes or edges should be more likely to be selected to host pressure sensors.  By ``central'' we mean a node or edge where water is more likely to pass through in order to meet clients' needs.  However,  to make this statement precise a more mathematical formulation is needed. 

In Complex Network Theory (CNT),  different concepts of centrality exist.  Examples of these are the betweennes,  closeness or degree centrality.  See \cite{newman2018networks} for an extensive overview on networks and centrality metrics.  A thorough discussion about all existing centralities would take us too far afield, but for the sake of clarity let us briefly introduce the betweennes centrality, a concept that in fact will be relevant for us in the following.

Given a node $i$ and a couple of other vertices $s$ and $t$,  say $m$ shortest paths connect $s$ to $t$. A fraction of these paths passes through $i$. The sum of all those fractions for every couple $(s, t)$ in the graph is called the betweennes centrality of the node $i$.  In more mathematical terms, we have the following definition:
\begin{equation}
C_i^B = \sum_{s \neq i \neq t \in \mathcal{V}} \frac{\sigma_{s,  t}(i)}{\sigma_{s,  t}} \, ,
\end{equation}
where $\sigma_{s,  t}$ is the number of shortest paths between $s$ and $t$, while $\sigma_{s,  t}(i)$ is the number of shortest paths through $i$. 

Remarkably,  an analogue formulation can be given for the centrality of an edge of the graph $G$.  In that case we speak of ``edge betweennes'' and, instead of considering the fraction of shortest paths passing through a node $i$, we compute the number of shortest paths linking two nodes $s$ and $t$ passing through an edge $e$ of the graph.  In particular, we take as a definition of the edge betweennes for an edge $e$:
\begin{equation}
C_e^B = \sum_{s \neq t \in \mathcal{V}} \frac{\sigma_{s,  t}(e)}{\sigma_{s,  t}} \, ,
\end{equation}
where $\sigma_{s,  t}$ is the number of shortest paths between $s$ and $t$ while $\sigma_{s,  t}(e)$ is the number of shortest paths between $s$ and $t$ passing through $e$.

A natural question is whether a definition like this would suffice to define centrality properties of a WDN. A possible answer to this question came from \cite{giustolisi2019tailoring} (see also \cite{giustolisi2017network}) where the authors tailored a suitable centrality metric for a generic WDN. Their analysis revolves essentially around three points \cite{giustolisi2019tailoring}, that we report here for clarity:
\begin{itemize}
\item edges (i.e.  pipes) are the most relevant components of a WDN, not nodes. Thus, centrality should be referred to edges rather than nodes,
\item pipes are characterized by asset features such as length, diameter and hydraulic resistance. Thus, a centrality should be weighted accordingly,
\item each node can represent a different component of the WDN: some of them are source nodes, others can be junctions (connections) or demand nodes. Thus, a centrality metric should be defined by taking into account where water flows start and end.
\end{itemize}

\begin{figure*}[!t]
\centering
\subfloat[Optimal sensor placement using the simulated annealing algorithm of the PyQUBO library. Red points mark nodes that have been selected to host a sensor. There are in total 48 red points and they appear to be uniformly distributed over the WDN with a tendency to aggregate in regions with a higher centrality.  35 out of the 48 sensors are found to be placed in correspondence of nodes classified as most accessible.  About 4.3$\%$ of the water flow passes through the selected nodes.]{\includegraphics[width=8.9cm]{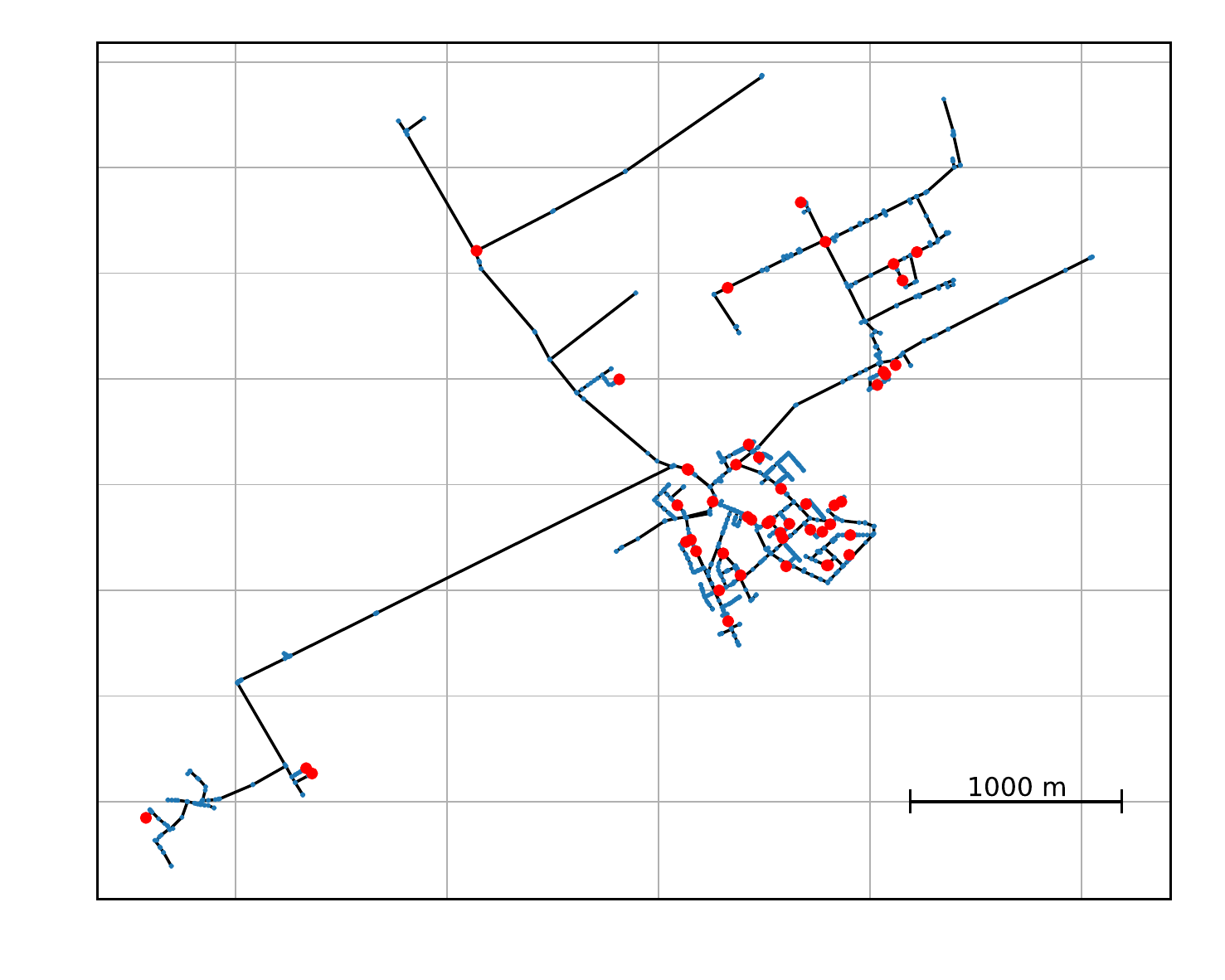}
\label{Optimal sensor placement using simulated annealing algorithms}}
\hfil
\subfloat[Energy spectral density of the simulation. The hyperparameters are chosen to be $(A, B, C, D) = (1, 30,5,1)$. On the horizontal axis we have the energies associated with different runs, while on the vertical axis the number of events normalized so to have total unit area under the histogram.  We have found that 100 runs is enough to get an acceptable solution to the sensors placement problem.]{\includegraphics[width=8.4cm]{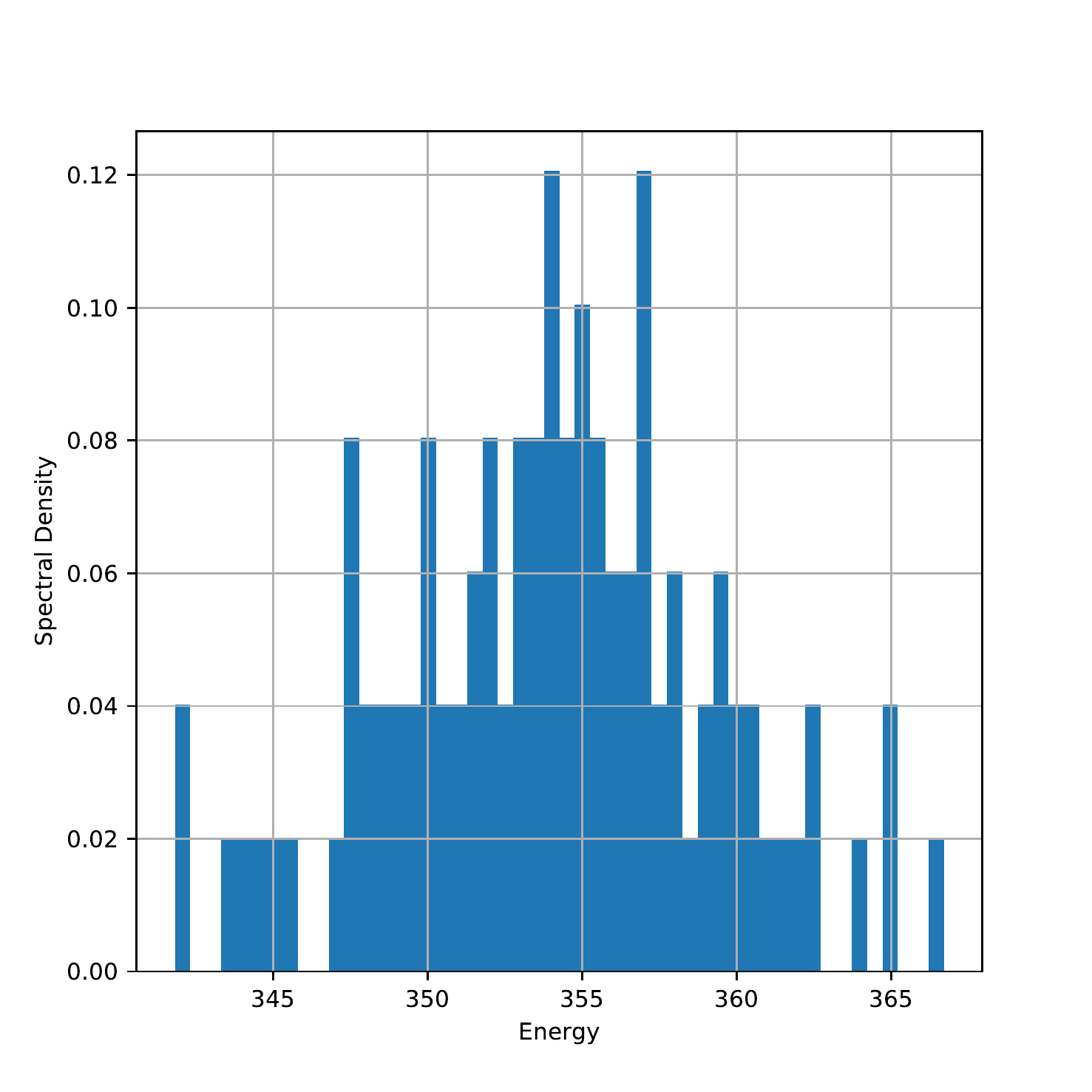}
\label{Energy spectral density of the simulation}}
\footnotesize{\caption{}}
\label{fig_sim}
\end{figure*}

While the first two items are relatively easy to deal with, the third requires some more explanation. What the authors of \cite{giustolisi2019tailoring} propose is to create a network of $n_f$ fictitious node around and connected only to source nodes,  with $n_f$ equal to\footnote{We are in fact oversimplifying a bit the analysis of \cite{giustolisi2019tailoring} where the authors also discuss the case where a tank undergoes an emptying/filling process, and directional devices -- installed near the tanks -- allow water to flow only in one direction.  The former case, in particular, is taken into account by replacing $n_f$ by $n_d V_T/V_D$, where $V_T$ and $V_D$ are the tank volume and average water volume supplied during  the operating cycle to demand nodes. In the present context,  we will not use this level of complexity.} (the integer part of) $ n/n_s$, where $n$ is the total number of nodes in the network and $n_s$ is the number of source nodes. In this way, source nodes behave like hubs for water supply and the centrality metric is weighted accordingly. 

In light of what we said, we propose to use as weights, $w_{ij}$,  the WDN tailored edge betweennes centrality as defined in \cite{giustolisi2019tailoring}. This will guarantee that the most central edges, in the sense just discussed, are more likely to be be assigned a sensor by the optimization algorithm.

\section{Case study: simulations and results on a real WDN}\label{Case study: simulations and results on a real WDN}

In this section,  we solve the sensor placement problem for the case of a real WDN. In order to minimize the Hamiltonian $H_P$ using classical and quantum annealing, we used PyQUBO, an open library for Python.  PyQUBO is an easy-to-use environment that allows to create and solve Ising or QUBO models and is fully integrated in D-Wave Ocean SDK. See \cite{zaman2021pyqubo} for a comprehensive overview on the subject. 

The WDN we use for the simulations is depicted in Fig.  \ref{wdn} and corresponds to a real water network in the Lombardy region in Italy.  It comprises $1368$ nodes (blue points),  one of which is associated with a tank (red point) and $1391$ edges (solid lines connecting nodes).  The WDN has a total length of approximately 26 kilometers, and serves about 4000 people.

The actual tank of our WDN lies in a very central position (red point in Fig. \ref{wdn}).  To make the problem more interesting, and for the purpose of this study only, we have added an extra outlying tank to the network (green point in Fig. \ref{wdn}). The resulting edge betweennes centrality is shown in Fig. \ref{wdn centrality}. 

The edge betweennes has been computed using pipe lengths as weights, and has been normalized to lie in the range 0 to 1.  As anticipated before, in order to take into account the position of the two tanks, a network of 683 fictitious nodes has been created around each tank.  Fig. \ref{wdn centrality} clearly shows that the most central edges (dark green) are those in proximity of the two tanks and where water is more likely to pass through to reach demand nodes, as we expected. A more refined version of the same centrality could be obtained by considering also other asset features, such as the pipe resistance or diameter.  

We have run our algorithm for many instances of hyperparameters $A$, $B$, $C$ and $D$.  Here, we report the case $(A, B, C, D) = (1, 30,5,1)$.  The vertex cost functions $f_i$ has been chosen to be $f_i = 1 - \widehat{v}_i$. Water consumption at each node $\widehat{v}_i$ corresponds to water need at each node over the period of one year.  We have also chosen to set $g_i = 0$ when the $i$-th node is considered easily accessible and $1$ otherwise.  Again, splitting the nodes in two categories of accessibility comes from engineering assessments.

Of course, more general choices of $f$ and $g$ exist, but for the purpose of the following analysis we will content ourselves with our choice.

The remaining of this section is made of two parts.  In the first part, we use simulated annealing methods to find the (approximate) global minimum of the Hamiltonian of the model of interest, eqn.  \eqref{final hamiltonian}, whereas in the second we use a hybrid quantum-classical approach on D-Wave Leap hybrid solver service (HSS). As it turns out, neither the Chimera nor the Pegasus D-Wave architectures were big enough to embed our instance problem  \eqref{final hamiltonian}. Thus, we drew on the D-Wave Hybrid Solver  where one should not worry about any embedding on D-Wave architectures nor should one think about choosing the right chain-strength. See the D-Wave documentation \cite{d-wave_doc} for more details on quantum and hybrid quantum-classical solvers.

Let us begin with the simulations performed using simulated annealing.

\subsection{Simulated Annealing}

As for the simulated annealing, we have used the Hamiltonian \eqref{final hamiltonian} with a number of sensors fixed at 48 units. The minimization is carried out performing 100 runs of the simulated annealing algorithm, for fixed hyperparameters $A$, $B$, $C$ and $D$, as described above.  Each run comes with an associated energy, each supposedly close to the global minimum, and the final sensor placement corresponds to the best result (minimum energy) among the 100 runs. The optimal sensor placement is shown in Fig. \ref{Optimal sensor placement using simulated annealing algorithms}.

Red points correspond to nodes selected by the optimization to host a sensor.  As we can see, sensors are distributed quite uniformly over the water network with a tendency to aggregate in regions with a higher centrality, ensuring a good spatial coverage. We have found that 35 out of 48 sensors are to be placed in correspondence of nodes classified as the most accessible, while the percentage of water flowing through the selected nodes is about $4.3 \%$ of the total flow.  The percentage of water that flows through the selected points might seem poor at first,  but in fact with a random sensor placement we would end up monitoring roughly the 3$\%$ of the total on average.  Thus, our configuration provides a consistent improvement. 

In Fig. \ref{Energy spectral density of the simulation} we show the energy spectral density of our minimization.  On the horizontal axis we give the energies associated with different runs, while on the vertical axis we depict the number of events normalized so to have total unit area under the histogram.  As we can see, to find the configuration closest to the global minimum, multiple runs are needed. We have found that 100 runs is enough to get an acceptable solution. The optimal solution is found when the energy is (approximately) $H_P = 341.8$.

As a final comment, we state that the minimization of the Hamiltonian in eqn. \eqref{final hamiltonian} for the case of the WDN shown in Fig. \ref{wdn} takes approximately $100$ seconds on a standard laptop equipped with a CPU Intel i7, which is indeed very satisfactory. In particular, this gives hope that running our algorithm on a bigger WDN could give sensible results in a reasonable time even on a classical computer.

Let us now move on to the case of the quantum-classical solver.

\subsection{Hybrid quantum-classical solver}
\begin{figure}[h!]
  \centering
   \includegraphics[width=9cm]{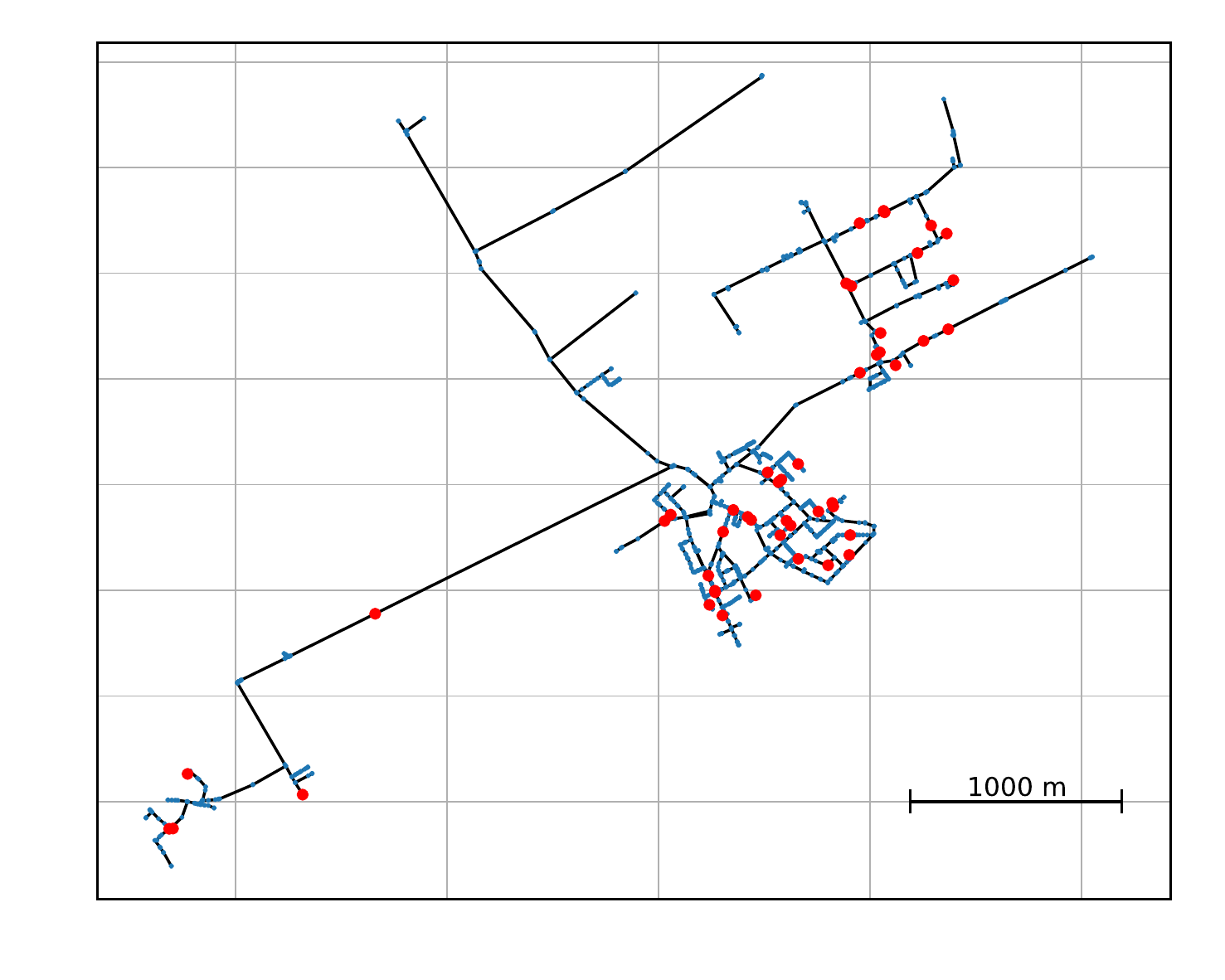}
   \caption{\footnotesize{Optimal sensor placement using the hybrid D-Wave Leap solver. Again, red points denote nodes that have been selected to host a sensor.  The number of selected nodes labeled as more accessible is of 39 out of 48 and percentage of water flowing through the selected notes is about $7.4 \%$ of the total flow. }}
   \label{wdn sensors placement hybrid 100}
\end{figure}

As for the hybrid quantum-classical annealing, we have run our algorithm on the D-Wave Leap hybrid solver service (HSS) \cite{d-wave_doc}.  The D-wave hybrid solver, differently from the quantum solver, does not need any embedding nor any choice of the chain strength \cite{d-wave_doc_2}.  The only tunable parameter is the Computation Time \cite{d-wave_doc_3}. We have found that the optimal solution for the minimization problem \eqref{final hamiltonian}, with $(A, B, C, D) = (1, 30,5,1)$, as above is achieved when the run-time is set at 100 seconds (similar to SA).

The outcome is slightly different from the case of the simulated annealing and is shown in Fig. \ref{wdn sensors placement hybrid 100}. The optimal solution is attained when the energy is approximately $H_P = 304.8$, lower than in the classical case.  We have also found that the number of selected nodes labeled as more accessible is of 39 out of 48 and percentage of water flowing through the selected notes is about $7.4 \%$ of the total flow. 

Even though we find, overall, that the hybrid solver performs better in this particular case, a systematic comparison between classical solvers and the hybrid quantum-classical approach to the optimization problem lies outside the main goals of the present paper. We hope to tackle these issues in future research.

\section{Conclusion}\label{Conclusion}

In this paper we considered a new quantum annealing heuristic method to solve the sensor placement problem on a generic Water Distribution Network (WDN).  Our method can be run on a classical computer by means of Simulated Annealing (SA) or on the publicly available D-Wave quantum solver.  Due to the limited number of qubits of the current D-Wave architectures, we relied on the hybrid quantum-classical solver reaching encouraging results.

The advantage of using an heuristic method that does not rely on hydraulic simulations lies in its flexibility and ease of use.  The algorithm is designed as a QUBO (or Ising problem) and is formulated in such a way sensors can be placed on a generic WDN as the result of a compromise between topological properties of the network and more distinguishing features, such as nodes accessibility or water consumption in a given area.

The hyper-parameters of the QUBO problem can be tuned so as to give more importance to either topology or clients' needs. We have given one instance of optimization for a set of hyper-parameters and compared results using simulated annealing and a hybrid quantum-classical method.

Our approach can be easily extended,  modified or generalized not only by adding or choosing different hyper-parameters for the QUBO problem, but also by considering different strategies for the sensor placement.  We hope that our work will inspire also how to tackle different optimization problems involving different spatial networks (transportation and mobility, electricity networks etc.). 

Moreover,  the proposed method proves to be promising to support, in real time, operational teams during the installation phase.  For example,  during the physical installation,  it is often necessary to relocate some of the sensors as a consequence of the fact that the positions indicated by the optimization algorithm may correspond to nodes that are difficult to access (lack of connectivity due to shielding,  lack of space,  etc...).  See Appendix \ref{A technical problem and its resolution} for a problem of this kind and its resolution.  Thus, the present work opens also the way to a more systematic assistance in the field thanks to the reduced computation time.

\section*{Acknowledgements}

It is a pleasure to thank David Amaro and Saskia Demulder for very interesting discussions and feedbacks.

\appendices

\section{Annealing Methods and Adiabatic Evolution}\label{Annealing Methods and Adiabatic Evolution}

In this appendix, we review some basic facts about annealing and adiabatic evolution.

“Annealing” is a term used in metallurgy and materials science.  It refers to the gradual cooling process of metal alloys and glassy materials to remove stress and defects.

The Simulated Annealing (SA) imitates annealing processes in computer simulations,  where a temperature is introduced into the optimization problem as thermal fluctuation \cite{kirkpatrick1983optimization}. The best solution can then be obtained by decreasing the temperature gradually.  Thermal excitations allow the system to escape local minima and reach states of lower and lower energy,  until it relaxes onto the ground state. 

Besides considering thermal fluctuations, we might consider quantum fluctuations as well.  An annealing process where quantum fluctuations drive the minimization is referrer to as Quantum Annealing (QA) \cite{kadowaki1998quantum, farhi2000quantum, kadowaki2002study}.  In the QA,  a quantum field plays a similar role to the temperature in the simulated annealing.  We start from a quantum-mechanical superposition of all possible states with equal weights, and then we evolve the system following the time-dependent Schroedinger equation. 

QA provides a formidable tool to find the lowest energy state encoding the solution to our combinatorial problem.  Excellent reviews on annealing methods and implementations can be found, for example, in \cite{nakahara2013lectures, tanaka2014physics}.

Let us say just a few words about the quantum adiabatic theorem \cite{messiah1999quantum, amin2009consistency}, a crucial ingredient in Adiabatic Quantum Optimization (AQC) and QA.

\subsection*{Quantum adiabatic theorem}

The time evolution of any quantum system is governed by the Schroedinger equation
\begin{equation}
i \hbar \frac{\partial \Psi}{\partial t} = H(X(t)) \Psi \, .
\end{equation}
Here, we assume that the Hamiltonian $H$ depends on a set of time-dependent parameters,  $X(t)$.  In order to discuss adiabatic changes of the quantum system,  the control parameters $X(t)$ must be varied slowly on the time-scale set by the energy-eigenvalues of the Hamiltonian $H$.  The set of control parameters $X$ can be thought of as parametrizing a topological space,  that we still denote as $X$.  For simplicity, $X$ is assumed to be connected, but not necessarily simply connected.  Examples of parameters that can be controlled externally are magnetic or electric fields, for instance. 

Here, $\Psi(t)$ denotes the state of the system at the time $t$ and,  mathematically,  is defined as a section of a bundle of Hilbert spaces over $X$, $\pi: \mathcal{H} \rightarrow X$.  Recall that a section of a bundle $\mathcal{H} \rightarrow X$ is a continuous map $\Psi: X \rightarrow \mathcal{H}$ such that $\pi \circ \Psi = \text{id}_X$, i.e.  $\pi(\Psi(x)) = x$, for any $x \in X$.  Another way to say this is that a section is a continuous assignment $x \mapsto \psi(x) \in \mathcal{H}_x$.  See \cite{nakahara2003geometry} for a review on fibre bundles and applications to Physics.

The inverse image, of the projection $\pi$,  $\pi^{-1}(x) = \mathcal{H}_x$ is the fibre at $x \in X$ and corresponds to a Hilbert space. In particular,  $\Psi(t) \in \mathcal{H}_{X(t)}$ is in a different fibre of the bundle for each $X(t)$. 

Of course, the bundle $\mathcal{H} \rightarrow X$ can be non-trivial. We assume that a trivialization is always possible. What this means is that over suitable open covers $\{ \mathcal{U}_{\alpha}\}$ of $X$ there are isomorphisms $\phi_{\alpha}: \pi^{-1}(\mathcal{U}_{\alpha}) \rightarrow \mathcal{U}_{\alpha} \times \mathcal{H}_0$, where $\mathcal{H}_0$ is some Hilbert space  which, in the finite dimensional case,  can be chosen to be isomorphic to $\mathbb{C}^{N}$,  $\mathcal{H}_0 \cong \mathbb{C}^{N}$, for some positive $N$ (e.g. $N = 2^n$ for a system of $n$ spins).  An important remark about trivializations is that if there are different trivializations of the bundle corresponding to different choices of basis $\{ \psi_{\alpha}(x)\}$, $\{ \psi_\beta(x)\}$, these should be related by some non-trivial map to the unitary group of $\mathcal{H}_0$, $U(\mathcal{H}_0)$, such that $\psi_{\beta} = U(\mathcal{H}_0)_{\beta \alpha}\psi_{\alpha}$.

The adiabatic theorem says,  in essence,  that a quantum system with slowly varying parameters follows, in the time evolution, the instantaneous eigenstate if we start the dynamics in one of the eigenstates of the initial Hamiltonian.  See \cite{messiah1999quantum, amin2009consistency} for full details.  Here, we are assuming that the degree of degeneracy of the energy spectrum does not change in the time evolution.  In particular,  if we start from the ground state (assumed to be unique for simplicity) of the initial Hamiltonian $H(0) = H_B$,  $\psi(0) = \psi_{g.s.}$,  then in the adiabatic approximation, the wave function stays close to the instantaneous ground state after a time $T$.  The adiabatic condition is often stated as
\begin{equation}\label{general adiabatic condition}
\text{max}_{0 \leq t \leq T} \frac{|\langle \psi_{n}(t),  \frac{\text{d}H(t)}{\text{d} t} \psi_{g.s.}(t) \rangle|}{E_{n, 0}(t)^2} \ll 1 \, , \qquad n \ge 1
\end{equation}
where $E_{n,0}(t) = E_n(t) - E_0(t)$,  while $\{ \psi_n\}$ is any trivialization of the Hilbert bundle at time $t$. Here, $n$ is a label for the $n$-th eigenstate (in order of increasing energy) of $H(t)$ with energy $E_n(t)$.  The bracket $\langle,  \rangle$ denotes an hermitian form on the fibre of the Hilbert bundle.

Equation \eqref{general adiabatic condition} can be rewritten for $n = 1$ (first excited state) as
\begin{equation}
T \gg \text{max}_{0 \leq u \leq 1} \frac{|\langle \psi_{1}(u),  \frac{\text{d}H(u)}{\text{d} u} \psi_{g.s.}(u) \rangle|}{E_{1, 0}(u)^2} \, ,
\end{equation}
which gives an estimate of what the total time evolution $T$ should be in order for the adiabatic condition to hold.

\section{A technical problem and its resolution}\label{A technical problem and its resolution}

Here, we would like to point out a technical problem that might arise when physically installing sensors on a WDN. 

Assume that our (classical or quantum) algorithm has found the optimal configuration for the sensor placement and $m (< s)$ sensors have been physically installed.  At the same time, assume that when installing the $(m+1)$th sensor a technical problem arises and the $(m+1)$th sensor cannot be installed at the selected position.  Running the algorithm again excluding that particular node would give in principle a completely different configuration and the work done for the first $m$ sensors should be undone. An unpleasant possibility.  Here we propose a simple way out for this problem. 

Define $\mathcal{V}_a \subset \mathcal{V}$, the subset of nodes where sensors have been physically installed and $\mathcal{V}_r \subset \mathcal{V}$ the subset of nodes that were selected by the algorithm but for one reason or another are not suitable to host a sensor. Of course,  $\mathcal{V}_a \cap \mathcal{V}_r = \emptyset$. 

We can run the minimization including a term of the form
\begin{equation}\label{good and bad points}
E \sum_{a \in \mathcal{V}_a} (x_a -1)^2 + E \sum_{r \in \mathcal{V}_r} x_r^2 \, ,
\end{equation}
where the first term is there to make sure that in the new optimization the nodes in $\mathcal{V}_a$ will be selected again (and no work has to be undone) while the second forbids the nodes in $\mathcal{V}_r$ to be selected in the new optimization.  $E$ is a positive real parameter.

We now define two vectors $k^{(a)}$, $k^{(r)} \in \mathbb{R}^n$, so that $k^{(a)}$ has non-zero entries (say 1's) in correspondence of nodes where a sensor has been installed (e.g.  if the 147th node hosts a sensor,  $k^{(a)}_{147} = 1$),  and $k^{(r)}$ has non-zero entries (again, 1's) in correspondence of nodes that should not be chosen by the optimization. It is an easy exercise to see that the term in \eqref{good and bad points} can be re-written as
\begin{equation}\label{good and bad points 2}
E \left[ k^{(a)} \cdot (k^{(a)} - x) + k^{(r)} \cdot x \right] \, ,
\end{equation}
where the dot is the Euclidean dot-product in $\mathbb{R}^n$ and the components of $x \in \mathbb{R}^n$ are the binary variables $x_i$ defined in the main text.

The advantage of considering eqn.  \eqref{good and bad points 2} rather than \eqref{good and bad points} is only for numerical purposes: if we initialize $k^{(a)}$, $k^{(r)}$ to be null vectors, we can just ``update'' them as explained above whenever a ``bad'' position (i.e. where sensors cannot be installed) is found on site and run again the optimization algorithm. 

\IEEEtriggeratref{28}
\bibliography{Bibliography}

\end{document}